\documentclass[11pt]{amsart2000} 
\usepackage{amssymb}
\newcommand{\LL}{{\it L}}
\newcommand{\I}{\mathcal{I}}
\newcommand{\A}{\mathcal{A}}
\newcommand{\B}{\mathcal{B}}
\newcommand{\C}{\mathcal{C}}
\newcommand{\D}{\mathcal{D}}
\newcommand{\V}{\mathcal{V}}
\newcommand{\U}{\mathcal{U}}
\newcommand{\X}{\mathcal{X}}
\newcommand{\Y}{\mathcal{Y}}
\newcommand{\OBC}{\mathcal{O}b(\mathcal{C})}
\newcommand{\OBD}{\mathcal{O}b(\mathcal{D})}

\newcommand{\MC}{\mathcal{M}(\mathcal{C})}
\newcommand{\MsC}{\mathcal{M}_\mathcal{C}}
\newcommand{\MD}{\mathcal{M}(\mathcal{D})}
\newcommand{\CCC}{ (\OBC, \omega, \MC, \mu, \circ)}
\newcommand{\DDD}{ (\OBD, \omega_\mathcal{D}, \MD, \mu_\mathcal{D},\circ)}
\newcommand{\CCCS}{ (\OBC, \omega', \MC, \mu', \circ)}
\newcommand{\tto}{\rightarrowtail}

\newcommand{\LFST}{{$\mathcal{F}SET(L)$}}
\newcommand{\LFSTT}{{$\mathcal{F}SET(L)^\top$}}
\newcommand{\LFTop}{{$\mathcal{F}TOP(L)$}}
\newcommand{\LFTOP}{{$\mathcal{F}TOP(L)^\top$}}
\newcommand{\impl}{\longmapsto}
\newtheorem{theorem}{Theorem}[subsection]

\newtheorem{definition}[theorem]{Definition}
\newtheorem{proposition}[theorem]{Proposition}
\newtheorem{corollary}[theorem]{Corollary}

\theoremstyle{definition}
\newtheorem{question}[theorem]{Question} 
\newtheorem{remark}[theorem]{Remark}
\newtheorem{example}[theorem]{Example}
\newtheorem{comments}[theorem]{Comments}
\numberwithin{equation}{theorem}

\begin{document} 
\setlength{\unitlength}{0.01in}
\linethickness{0.01in}
\begin{center}
\begin{picture}(474,66)(0,0) 
\multiput(0,66)(1,0){40}{\line(0,-1){24}}
\multiput(43,65)(1,-1){24}{\line(0,-1){40}}
\multiput(1,39)(1,-1){40}{\line(1,0){24}}
\multiput(70,2)(1,1){24}{\line(0,1){40}}
\multiput(72,0)(1,1){24}{\line(1,0){40}}
\multiput(97,66)(1,0){40}{\line(0,-1){40}} 
\put(143,66){\makebox(0,0)[tl]{\footnotesize Proceedings of the Ninth Prague Topological Symposium}}
\put(143,50){\makebox(0,0)[tl]{\footnotesize Contributed papers from the symposium held in}}
\put(143,34){\makebox(0,0)[tl]{\footnotesize Prague, Czech Republic, August 19--25, 2001}}
\end{picture}
\end{center}
\vspace{0.25in}
\setcounter{page}{271}
\title[Fuzzy functions and $L$-Top]{Fuzzy functions and an extension of
the category $L$-Top of Chang-Goguen $L$-topological spaces}
\author{Alexander P. \v{S}ostak} 
\address{Department of Mathematics\\
University of Latvia\\
Riga\\
Latvia}
\email{sostaks@com.latnet.lv}
\thanks{Partly supported by grant 01.0530 of Latvijas Zin\=atnes Padome}
\thanks{This article will be revised and submitted for publication
elsewhere.}
\thanks{Alexander P. \v{S}ostak,
{\em Fuzzy functions and an extension of the category $L$-Top of 
Chang-Goguen $L$-topological spaces},
Proceedings of the Ninth Prague Topological Symposium, (Prague, 2001),
pp.~271--294, Topology Atlas, Toronto, 2002}
\begin{abstract}
We study $\mathcal{F}TOP(L)$, a fuzzy category 
with fuzzy functions in the role of morphisms. This category has the
same objects as the category L-TOP of Chang-Goguen L-topological spaces,
but an essentially wider class of morphisms---so called fuzzy functions
introduced earlier in our joint work with U. H\"ohle and H. Porst.
\end{abstract}
\maketitle 

\section*{Introduction} 

In research works where fuzzy sets are involved, in particular, in 
Fuzzy Topology, mostly certain usual functions are taken as 
morphisms: they can be certain mappings between corresponding sets, or
between the fuzzy powersets of these sets, etc. 
On the other hand, in our joint works with U.~H\"ohle and H.E.~Porst
\cite{HPS1}, \cite{HPS2} a certain class of $L$-relations (i.e.\ mappings
$F: X\times Y \to L$) was distinguished which we view as ($L$-){\it fuzzy
functions} from a set $X$ to a set $Y$; these fuzzy functions play the
role of morphisms in an {\it $L$-fuzzy category} of sets \LFST,
introduced in \cite{HPS2}.

Later on we constructed a fuzzy category \LFTop\ related to topology
with fuzzy functions in the role of morphisms, see \cite{So2000}. 
Further, in \cite{So2001} a certain uniform counterpart of \LFTop\ was
introduced.
Our aim here is to continue the study of \LFTop. 
In particular, we show that the top frame \LFTop$^\top$ of the fuzzy
category \LFTop\ is a topological category in H. Herrlich's sense
\cite{AHS}) over the top frame \LFST$^\top$ of the fuzzy category \LFST.

In order to make exposition self-contained, we start with Section 
1 Prerequisities, where we briefly recall the three basic concepts which
are essentially used in this work: they are the concepts of a $GL$-monoid
(see e.g.\ \cite{Ho91}, \cite{Ho94}, etc.); 
of an $L$-valued set (see e.g.\ \cite{Ho92}, etc.), and of an 
$L$-fuzzy
category (see e.g.\ \cite{So91}, \cite{So92}, \cite{So97}, etc.). 
In Section 2 we consider basic facts about fuzzy functions and introduce
the $L$-fuzzy category \LFST\ \cite {HPS1}, \cite{HPS2}. 
The properties of this fuzzy category and some related categories are the
subject of Section 3.
\LFST\ is used as the ground category for the $L$-fuzzy category \LFTop\ 
whose objects are Chang-Goguen $L$-topological spaces \cite{Ch},
\cite{Go73}, and whose morphisms are certain fuzzy functions,
i.e.\ morphisms from \LFST. Fuzzy category \LFTop\ is considered in
Section 4. 
Its crisp top frame \LFTOP\ is studied in Section 5. 
In particular, it is shown that \LFTOP\ is a topological category over
\LFST$^\top$. 
Finally, in Section 6 we consider the behaviour of the topological
property of compactness with respect to fuzzy functions --- in other
words in the context of the fuzzy category \LFTop\ and, specifically, in
the context of the category \LFTOP.

\subsection*{Acknowledgements} 
The author is thankful to U.~H\"ohle for many useful discussions which
stimulated this work.

\section{Prerequisities} 

\subsection{$GL$-monoids}

Let $(L, \leq)$ be a complete infinitely distributive lattice, i.e.\ 
$(L, \leq)$ is a partially ordered set such that for every subset $A
\subset L$ the join $\bigvee A$ and the meet $\bigwedge A$ are defined
and 
$(\bigvee A) \wedge \alpha = \bigvee \{ a\wedge \alpha) \mid a \in A \}$
and 
$(\bigwedge A) \vee \alpha = \bigwedge \{a\vee \alpha) \mid a \in A \}$
for every $\alpha \in L$.
In particular, 
$\bigvee L =: \top$ and $\bigwedge L =: \bot$ 
are respectively the universal upper and the universal lower bounds 
in $L$. 
We assume that $\bot \ne \top$, i.e.\ $L$ has at least two elements.

A $GL-$monoid (see \cite{Ho91}, \cite{Ho92}, \cite{Ho94}) is a complete 
lattice enriched with a further binary operation 
$*$, i.e.\ a triple $(L, \leq, *)$ such that:
\begin{enumerate}
\item[(1)] 
$*$ is monotone, i.e.\ $\alpha \leq \beta$ implies 
$\alpha * \gamma \leq \beta * \gamma$, 
$\forall \alpha, \beta, \gamma \in \LL$;
\item[(2)] 
$*$ is commutative, i.e.\ $\alpha * \beta = \beta * \alpha$, 
$\forall \alpha, \beta \in \LL$;
\item[(3)] 
$*$ is associative, i.e.\ 
$\alpha * (\beta * \gamma) = (\alpha * \beta) * \gamma$,
$\forall \alpha, \beta, \gamma \in L$;
\item[(4)] 
$(L,\leq,*)$ is integral, i.e.\ $\top$ acts as the unity: 
$\alpha * \top = \alpha$, $\forall \alpha \in \LL$;
\item[(5)] 
$\bot$ acts as the zero element in $(L, \leq, *)$, 
i.e.\ $\alpha * \bot = \bot$, $\forall \alpha \in \LL$;
\item[(6)] 
$*$ is distributive over arbitrary joins, i.e.\
$\alpha * (\bigvee_j \beta_j) = \bigvee_j (\alpha * \beta_j)$, 
$\forall \alpha \in \LL, \forall \{ \beta_j : j \in J \} \subset \LL$;
\item[(7)] 
$(L, \leq, *)$ is divisible, i.e.\ $\alpha \leq \beta$ implies 
existence of $\gamma \in L$ such that $\alpha = \beta * \gamma$.
\end{enumerate}

It is known that every $GL-$monoid is residuated, i.e.\ there exists a 
further binary operation ``$\impl$'' (implication) on $L$ satisfying the 
following condition:
$$\alpha * \beta \leq \gamma \Longleftrightarrow 
\alpha \leq (\beta \impl \gamma) 
\qquad \forall \alpha, \beta, \gamma \in L.$$
Explicitly implication is given by
$$\alpha \impl \beta = 
\bigvee \{ \lambda \in L \mid \alpha * \lambda \leq \beta \}.$$

Below we list some useful properties of $GL-$monoids 
(see e.g.\ \cite{Ho91}, \cite{Ho92}, \cite{Ho94}):
\begin{enumerate}
\item[(i)] 
$\alpha \impl \beta = \top \Longleftrightarrow \alpha \leq \beta$;
\item[(ii)] 
$\alpha \impl (\bigwedge_i \beta_i) = \bigwedge_i (\alpha \impl \beta_i)$;
\item[(iii)] 
$(\bigvee_i \alpha_i) \impl \beta = \bigwedge_i (\alpha_i \impl \beta)$;
\item[(v)] 
$\alpha * (\bigwedge_i \beta_i) = \bigwedge_i (\alpha * \beta_i)$;
\item[(vi)] 
$(\alpha \impl \gamma ) * (\gamma \impl \beta) \leq \alpha \impl \beta$;
\item[(vii)] 
$\alpha * \beta \leq (\alpha * \alpha) \vee (\beta * \beta)$.
\end{enumerate}

Important examples of $GL$-monoids are Heyting algebras and 
$MV$-alg\-ebras. 
Namely, a {\it Heyting algebra} is $GL$-monoid of the type 
$(L,\leq,\wedge,\vee,\wedge)$ (i.e.\ in case of a Heyting algebra 
$\wedge = *$), cf.\ e.g.\ \cite{Jhst}. 
A $GL$-monoid is called an {\it $MV$-algebra} if 
$(\alpha \impl \bot) \impl \bot = \alpha \quad \forall \alpha \in L$, 
\cite{Ch58}, \cite{Ch59}, see also \cite[Lemma 2.14]{Ho94}. 
Thus in an $MV$-algebra an order reversing involution $^c: L \to L$ can
be naturally defined by setting 
$\alpha^c := \alpha \impl \bot \quad \forall \alpha \in L$.

If $X$ is a set and $L$ is a $GL$-monoid, then the fuzzy powerset 
$L^X$ in an obvious way can be pointwise endowed with a structure 
of a $GL$-monoid. 
In particular the $L$-sets $1_X$ and $0_X$ defined by $1_X (x):= \top$
and $0_X (x) := \bot$ $\forall x \in X$ are respectively the universal
upper and lower bounds in $L^X$.

In the sequel $L$ denotes an arbitrary $GL$-monoid.

\subsection{$L$-valued sets}

Following U.~H\"ohle (cf.\ e.g.\ \cite{Ho92}) by a (global) {\it 
$L$-valued set} we call a pair $(X,E)$ where $X$ is a set and $E$ 
is an {\it $L$-valued equality}, i.e.\ a mapping $E: X \times X \to 
L$ such that
\begin{enumerate} 
\item[(1eq)] 
$E(x,x) = \top$
\item[(2eq)]
$E(x,y) = E(y,x) \quad \forall x, y \in X$;
\item[(3eq)] 
$E(x,y)*E(y,z) \leq E(x,z) \quad \forall x, y, z \in X$.
\end{enumerate}

A mapping $f: (X,E_X) \to (Y,E_Y)$ is called {\it extensional} if 
$$E_X(x,x') \leq E_Y(f(x),f(x'))\ \forall x,x' \in X.$$
Let $SET(L)$ denote the category whose objects are $L$-valued sets and
whose morphisms are extensional mappings between the corresponding
$L$-valued sets.
Further, recall that an $L$-set, or more precisely, an $L$-subset of a 
set $X$ is just a mapping $A: X \to L$.
In case $(X,E)$ is an $L$-valued set, its $L$-subset $A$ is called {\it
extensional} if 
$$\bigvee_{x\in X} A(x) * E(x,x') \leq A(x') \quad \forall x' \in X.$$

\subsection{\LL-fuzzy categories} 

\begin{definition}[\cite{So91}, \cite{So92}, \cite{So97}, \cite{So99}]
An \LL-fuzzy category is a quintuple $\C = \CCC$ where 
$\C_\bot = (\OBC, \MC, \circ)$ is a usual (classical) category 
called {\em the bottom frame} of the fuzzy category $\C$; 
$\omega: \OBC \longrightarrow \LL $ is an $L$-subclass of the class of 
objects $\OBC$ of $\C_\bot$ and $\mu : \MC \longrightarrow \LL $ is an
$L$-subclass of the class of morphisms $\MC$ of $\C_\bot$. 
Besides $\omega$ and $\mu$ must satisfy the following conditions:
\begin{enumerate}
\item[(1)] 
if $f: X \to Y$, then $ \mu (f) \leq \omega (X) \wedge \omega (Y)$;
\item[(2)] 
$\mu (g \circ f) \geq \mu (g) * \mu (f)$ whenever composition $g \circ f$
is defined;
\item[(3)] 
if $e_X: X \to X$ is the identity morphism, then $\mu(e_X) = \omega(X)$.
\end{enumerate}
\end{definition}

Given an \LL-fuzzy category $ \C = \CCC$ and $X \in \OBC$, the intuitive
meaning of the value $\omega (X)$ is the {\it degree} to which a
potential object $X$ of the \LL-fuzzy category $\C$ is indeed its 
object; similarly, for $f \in \MC$ the intuitive meaning of $\mu (f)$ is
the degree to which a potential morphism $f$ of $\C$ is indeed its
morphism.

\begin{definition} 
Let $\C = \CCC$ be an \LL-fuzzy category. 
By an (\LL-fuzzy) subcategory of $\C$ we call an \LL-fuzzy category 
$$\C' = \CCCS$$ 
where $\omega' \leq \omega$ and $\mu' \leq \mu$.
A subcategory $\C'$ of the category $\C$ is called full if 
$\mu'(f) = \mu(f) \wedge \omega'(X) \wedge \omega'(Y)$ for every 
$f \in \MsC (X,Y)$, and all $X$,$Y \in \OBC$.
\end{definition}

Thus an \LL-fuzzy category and its subcategory have the same classes of 
potential objects and morphisms. 
The only difference of a subcategory from the whole category is in
\LL-fuzzy classes of objects and morphisms, i.e.\ in the belongness 
degrees of potential objects and morphisms. 

Let $\C = (\OBC, \MC, \circ)$ be a crisp category and 
$\D = (\OBD, \MD, \circ)$ be its subcategory. 
Then for every $GL$-monoid \LL\ the category $\D$ can be identified with
the \LL-fuzzy subcategory 
$$\tilde{\D} = (\OBC, \omega', \MC, \mu', \circ)$$ 
of $\C$ such that 
$\omega'(X) = \top$ if $X \in \OBD$ and 
$\omega'(X) = \bot$ otherwise; 
$\mu' (f) = \top$ if $f \in \MD$ and $\mu' (f) = \bot$ otherwise. 
In particular, ${\tilde{\D}}_\top = \D$.

On the other hand sometimes it is convenient to identify a fuzzy
subcategory 
$$\C' = \CCCS$$
of the fuzzy category 
$$\C = \CCC$$ 
with the fuzzy category 
$$\D = \DDD$$
where 
$$\OBD := \{X \in \OBC \mid \omega'(X) \ne \bot \}$$ 
and 
$$\MD := \{f \in \MC \mid \mu'(f) \ne \bot \}$$ 
and $\omega_\mathcal{D}$ and $\mu_\mathcal{D}$ are restrictions of
$\omega'$
and $\mu'$ to $\OBD$ and $\MD$ respectively.

\section{Fuzzy functions and fuzzy category \LFST} 

As it was already mentioned above, the concept of a fuzzy function and 
the corresponding fuzzy category \LFST\ were introduced in \cite{HPS1}, 
\cite{HPS2}. 
There were studied also basic properties of fuzzy functions. 
In this section we recall those definitions and results from \cite{HPS2}
which will be needed in the sequel\footnote{Actually, the subject of
\cite{HPS1}, \cite{HPS2} was a more general fuzzy category $L$-\LFST\
containing \LFST\ as a full subcategory.
However, since for the merits of this work the category \LFST\ is of 
importance, when discussing results from \cite{HPS1}, \cite{HPS2} we
reformulate (simplify) them for the case of \LFST\ without mentioning
this every time explicitly}. 
Besides some new needed facts about fuzzy functions will be established
here, too. 

\subsection{Fuzzy functions and category \LFST}

\begin{definition}[cf.\ {\cite[2.1]{HPS1}}]
A fuzzy function\footnote{Probably, the name {\it an $L$-fuzzy 
function} would be more adequate here. 
However, since the $GL$-monoid $L$ is considered to be fixed, and since
the prefix ``$L$'' appears in the text very often, we prefer to say just
{\it a fuzzy function}} 
$F$ from an $L$-valued set $(X,E_X)$ to $(Y,E_Y)$ (in symbols
$F: (X,E_X) \rightarrowtail (Y,E_Y))$) is a mapping 
$F: X \times Y \to L$ such that 
\begin{enumerate}
\item[(1ff)] 
$F(x,y) * E_Y(y,y') \leq F(x,y') 
\quad \forall x \in X, \forall y,y' \in Y$;
\item[(2ff)] 
$E_X(x,x') * F(x,y) \leq F(x',y) 
\quad \forall x,x' \in X, \forall y \in Y$; 
\item[(3ff)] 
$F(x,y) * F(x,y') \leq E_Y(y,y') 
\quad \forall x \in X, \forall y,y' \in Y$.
\end{enumerate}
\end{definition}

{\small Notice that conditions (1ff)--(2ff) say that $F$ is a 
certain $L-$relation, while axiom (3ff) together with evaluation 
$\mu(F)$ (see Subsection \ref{fsetl}) specify that the $L$-relation 
$F$ is a fuzzy {\it function}}.

\begin{remark}
Let $F: (X,E_X) \tto (Y,E_Y)$ be a fuzzy function, $X' \subset X$, 
$Y' \subset Y$, and let the $L$-valued equalities $E_{X'}$ and $E_{Y'}$
on $X'$ and $Y'$ be defined as the restrictions of the equalities $E_X$
and $E_Y$ respectively. 
Then defining a mapping $F': X'\times Y' \to L$ by the equality
$F'(x,y) = F(x,y)\ \forall x \in X', \forall y \in Y'$ 
a fuzzy function $F': (X',E_{X'}) \tto (Y',E_{Y'})$ is obtained. 
We refer to it as the {\it restriction of $F$} to the subspaces 
$(X',E_{X'})$ $(Y',E_{Y'})$
\end{remark}

Given two fuzzy functions $F: (X,E_X) \tto (Y,E_Y)$ and 
$G: (Y,E_Y) \tto (Z,E_Z)$ we define their {\it composition} 
$G \circ F: (X,E_X) \tto (Z,E_Z)$ by the formula
$$(G\circ F)(x,z) = \bigvee_{y\in Y} \bigl( F(x,y) * G(y,z))\bigr).$$ 

In \cite{HPS2} it was shown that the composition $G \circ F$ is indeed a 
fuzzy function and that the operation of composition is associative.
Further, if we define the identity morphism by the corresponding 
$L$-valued equality:
$E_X: (X,E_X) \tto (X,E_X),$ we come to a category \LFST\ 
whose objects are $L$-valued sets and whose morphisms are fuzzy functions
$F: (X,E_X) \tto (Y,E_Y)$. 
 
\subsection{Fuzzy category \LFST}\label{fsetl}

Given a fuzzy function $F: (X,E_X) \tto (Y,E_Y)$ let
$$\mu(F) = \inf_x \sup_y F(x,y).$$
Thus we define an $L$-subclass $\mu$ of the class of all morphisms of
\LFST. 
In case $\mu(F) \geq \alpha$ we refer to $F$ as a {\it fuzzy 
$\alpha$-function}.

If $F: (X,E_X) \tto (Y,E_Y) \mbox{ and } G: (Y,E_Y) \tto (Z,E_Z)$ are 
fuzzy functions, then $\mu (G \circ F) \geq \mu(G) * \mu(F)$ \cite{HPS2}.
Further, given an $L$-valued set $(X,E)$ let 
$\omega(X,E) := \mu(E) = \inf_x E(x,x) = \top$.
Thus a {\it fuzzy category} \LFST = $(FSET(L), \omega, \mu)$ is
obtained.

\begin{example}\label{ex-crisp}
Let $* = \wedge$ and $E_Y$ be a crisp equality on $Y$, i.e.\ 
$E_Y(y,y') = \top \mbox{ iff } y = y'$, and 
$E_Y(y,y') = \bot$ otherwise. 
Then every fuzzy function $F: (X,E_X) \tto (Y,E_Y)$ such that 
$\mu(F) = \top$ is uniquely determined by a usual function 
$f: X \to Y$. 
Indeed, let $f(x) = y \mbox{ iff } F(x,y) = \top$. 
Then condition (3ff) implies that there cannot be $f(x)=y$, $f(x) = y'$
for two different $y$, $y' \in Y$ and condition $\mu(F) = \top$ 
guarantees that for every $x \in X$ one can find $y \in Y$ such 
that $f(x)=y$.
If besides $E_X$ is crisp, then, vice versa, every mapping $f: X \to Y$
can be viewed as a fuzzy mapping $F: (X,E_X) \tto (Y,E_Y)$ (since the
conditions of extensionality (2ff) and (3ff) are automatically fulfilled
in this case)
\end{example}

\begin{remark}
If $F'\colon (X',E_{X'}) \tto (Y,E_{Y})$ is the restriction of 
$F\colon (X,E_X) \tto (Y,E_Y)$ (see Remark above) and 
$\mu(F) \geq \alpha$, then $\mu(F') \geq \alpha$. 
However, generally the restriction 
$F'\colon (X',E_{X'}) \tto (Y',E_{Y'})$ of $F\colon (X,E_X) \tto (Y,E_Y)$
may fail to satisfy condition $\mu(F') \geq \alpha$.
\end{remark}

\subsection{Images and preimages of $L$-sets under fuzzy
functions}\label{impref}

Given a fuzzy function $F: (X, E_X) \tto (Y, E_Y)$ and $L$-subsets
$A: X \to L$ and $B: Y \to L$ of $X$ and $Y$ respectively, we define the
fuzzy set $F^{\rightarrow}(A): Y \to L$ (the image of $A$ under $F$) by
the equality $F^{\rightarrow}(A)(y) = \bigvee_x F(x,y) * A(x)$ and the
fuzzy set $F^{\leftarrow}(B): X \to L$ (the preimage of $B$ under $F$) 
by the equality $F^{\leftarrow}(B)(x) = \bigvee_y F(x,y) * B(y)$.

Note that if $A \in L^X$ is extensional, then 
$F^{\rightarrow}(A) \in L^Y$ is extensional (by (2ff)) and if 
$B \in L^Y$ is extensional, then $F^{\leftarrow}(B) \in L^X$ is 
extensional (by (3ff)).

\begin{proposition}[Basic properties of images and preimages of $L$-sets
under fuzzy functions]
\label{im-pr}
\mbox{}
\begin{enumerate}
\item 
$F^{\rightarrow}(\bigvee_{i \in \mathcal{I}} (A_i) = 
\bigvee_{i \in \mathcal{I}} F^{\rightarrow}(A_i) 
\qquad \forall \{A_i: i \in {\mathcal{I}} \} \subset L^X$;
\item 
$F^{\rightarrow}(A_1\bigwedge A_2) \leq 
F^{\rightarrow}(A_1) \bigwedge F^{\rightarrow}(A_2) 
\qquad 
\forall A_1, A_2 \in L^X$;
\item 
If $L$-sets $B_i$ are {\em extensional}, then 
$$
\bigwedge_{i \in {\mathcal{I}}} 
F^{\leftarrow}(B_i)*\mu(F)^2 
\leq 
F^{\leftarrow}(\bigwedge_{i \in \mathcal{I}} B_i) 
\leq 
\bigwedge_{i \in \mathcal{I}} F^{\leftarrow}(B_i) 
\qquad \forall \{B_i: i \in \mathcal{I} \} \subset L^Y.
$$
In particular, if $\mu(F) = \top$, then 
$F^{\leftarrow}(\bigwedge_{i\in \mathcal{I}} B_i) = 
\bigwedge_{i \in \mathcal{I}} F^{\leftarrow}(B_i)$ 
for every family of extensinal $L$-sets 
$\{B_i: i \in \mathcal{I} \} \subset L^Y$. 
\item 
$F^{\leftarrow}(\bigvee_{i \in \mathcal{I}} B_i) = 
\bigvee_{i \in \mathcal{I}} F^{\leftarrow}(B_i) 
\qquad \forall \{B_i: i \in {\mathcal{I}} \} \subset L^Y$.
\item 
$A*\mu(F)^2 \leq F^{\leftarrow}(F^{\rightarrow}(A))$ 
for every $A \in L^X$.
\item 
$F^{\to}\bigl(F^{\gets}(B)\bigr) \leq B$ 
for every {\em extensional} $L$-set $B\in L^Y$.
\item 
$F^{\leftarrow}(c_Y) \geq \mu(F) * c$ where $c_Y: Y \to L$ is the constant
function taking value $c \in L$. 
In particular, $F^{\leftarrow}(c_Y) = c$ if $\mu(F) = \top$.
\end{enumerate}
\end{proposition}

\begin{proof}
(1).
$$
\begin{array}{lll}
\bigl(\bigvee_i F^{\rightarrow}(A_i)\bigr)(y)&
=&
\bigvee_i \bigvee_x \bigl(F(x,y) * A_i(x) \bigr)\\
&
=&
\bigvee_x \bigvee_i \bigl(F(x,y) * A_i (x)\bigr)\\
&
=&
\bigvee_x (F(x,y) * (\bigvee_i A_i)(x))\\
&
=&
F^{\rightarrow}(\bigvee_i A_i)(y).
\end{array}
$$

(2).
The validity of (2) follows from the monotonicity of $F^{\to}$.

(3).
To prove property 3 we first establish the following inequality 
\begin{equation}\label{basic1}
\bigvee_{y\in Y}\bigl(F(x,y)\bigr)^2 \geq 
\bigl(\bigvee_{y\in Y} F(x,y)\bigr)^2.
\end{equation}
Indeed, by a property (vii) of a GL-monoid 
$$
\begin{array}{lll}
\Bigl(\bigvee_{y\in Y} F(x,y)\Bigr)^2&
=&
\Bigl(\bigvee_{y\in Y} F(x,y)\Bigr) * \Bigl(\bigvee_{y'\in Y}
F(x,y')\Bigr)\\
&
=&
\bigvee_{y,y'\in Y} \Bigl(F(x,y) * F(x,y') \Bigr)\\
&
\leq&
\bigvee_{y,y'\in Y} \Bigl(F(x,y)^2 \vee F(x,y')^2\Bigr)\\
&
=&
\bigvee_{y\in Y} \Bigl(F(x,y)\Bigr)^2.
\end{array}$$
In particular, it follows from (\ref{basic1}) that 
\begin{equation}\label{basic2}
\forall x \in X\ 
\bigvee_{y\in Y} \bigl(F(x,y)\bigr)^2 \geq \mu(F)^2.
\end{equation}
Now, applying (\ref{basic2}) and taking into account extensionality
of
$L$-sets $B_i$, we proceed as follows:
$$
\begin{array}{lll}
\multicolumn{3}{l}{
\Bigl(\bigwedge_i F^{\gets}(B_i)\Bigr)(x) * \bigl(\mu(F)\bigr)^2
}
\\
&
\leq&
\Bigl(\bigwedge_i \bigl(\bigvee_{y_i \in Y} (F(x,y_i)*B_i(y_i)\bigr)\Bigr) * 
\bigvee_{y\in Y}\bigl(F(x,y)\bigr)^2\\
&
=&
\bigvee_{y\in Y} \Bigl(\bigl(F(x,y)\bigr)^2 * 
\bigwedge_i\bigl(\bigvee_{y_i \in Y} F(x,y_i) * B_i(y_i)\bigr)\\
&
=&
\bigvee_{y\in Y}\Bigl(F(x,y) * \bigl(\bigwedge_i \bigl(F(x,y) * 
\bigl(\bigvee_{y_i \in Y} F(x,y_i) * B_i(y_i)\bigr)\bigr)\Bigr)\\
&
=&
\bigvee_{y\in Y}\Bigl(F(x,y) * \bigl(\bigwedge_i \bigl(\bigvee_{y_i\in 
Y} (F(x,y) * F(x,y_i)\bigr) \bigr) * B_i(y_i)\bigr)\Bigr)\\
&
\leq&
\bigvee_{y\in Y}\Bigl(F(x,y) * \bigl(\bigwedge_i \bigl(\bigvee_{y_i \in Y} 
E(y,y_i) * B_i(y_i)\bigr)\bigr)\Bigr)\\
&
\leq&
\bigvee_{y\in Y}\bigl(F(x,y) * \bigl(\bigwedge_i B_i(y)\bigr)\bigr)\\
&
=&
F^{\gets}\bigl(\bigwedge_i B_i)\bigr)(x),$$
\end{array}
$$
and hence 
$$\Bigl(\bigwedge_i F^{\gets}(B_i)\Bigr)* \bigl(\mu(F)\bigr)^2 
\leq F^{\gets}\bigl(\bigwedge_i B_i)\bigr).$$ 
To complete the proof notice that the inequality
$$F^{\leftarrow}(\bigwedge_{i \in \mathcal{I}} B_i) \leq 
\bigwedge_{i \in \mathcal{I}} F^{\leftarrow}(B_i)$$ 
is obvious. 

(4).
The proof of (4) is similar to the proof of (1) and is therefore 
omitted.

(5).
Let $A \in L^X$, then 
$$
\begin{array}{lll}
F^{\leftarrow}(F^{\rightarrow}(A))(x)&
=&
\bigvee_y (F(x,y)*F^{\rightarrow}(A)(y))\\
&
=&
\bigvee_y \bigl(F(x,y)*\bigl(\bigvee_{x'} F(x',y) * A(x') \bigr)\bigr)\\
&
\geq&
\bigvee_y (F(x,y)^2 * A(x))\\
&
\geq&
(\mu(F))^2 * A(x) 
\end{array}
$$
for every $x\in X$, and hence 5 holds.

(6).
To show property 6 assume that $B\in L^Y$ is extensional. 
Then
$$
\begin{array}{lll}
F^{\rightarrow}(F^{\leftarrow}(B))(y)&
=&
\bigvee_x\bigl(F(x,y) * F^{\gets}(B)(x)\bigr)\\
&
=&
\bigvee_x F(x,y) * \bigl(\bigvee_{y'} F(x,y') * B(y') \bigr)\\
&
=&
\bigvee_{x \in X, y\in Y} \bigl(F(x,y) * F(x,y') * B(y')\bigr)\\
&
\leq&
E_Y(y,y')*B(y')\\
&
\leq&
B(y), 
\end{array}
$$
and hence 
$F^{\to}\bigl(F^{\gets}(B)\bigr) \leq B$.

(7).
The proof of property 7 is straightforward and therefore omitted.
\end{proof}

\begin{comments}
\mbox{}
\begin{enumerate}
\item 
Properties 1,2 and 4 were proved in \cite[Proposition 3.2]{HPS2}. 
Here we reproduce these proofs in order to make the article 
self-contained. 
\item 
The inequality in item 2 of the previous proposition obviously 
cannot be improved even in the crisp case. 
\item 
One can show that the condition of extensionality cannot be omitted 
in items 3 and 6. 
\item 
The idea of the proof of Property 3 was communicated to the author by
U.~H\"ohle in Prague at TOPOSYM in August 2001. 
\item 
In \cite{HPS2} there was established the following version of 
Property 3 without the assumption of extensionality of $L$-sets 
$B_i$ in case $L$ is completely distributive:
$$ {(\bigwedge_{i \in \mathcal{I}} F^{\leftarrow}(B_i))}^5 \leq 
F^{\leftarrow}(\bigwedge_{i \in \mathcal{I}} B_i) \leq 
\bigwedge_{i \in \mathcal{I}} F^{\leftarrow}(B_i) \qquad 
\forall \{B_i: i \in \mathcal{I} \} \subset L^Y$$
and
$$ \bigwedge_{i \in \mathcal{I}} F^{\leftarrow}(B_i) = 
F^{\leftarrow}(\bigwedge_{i \in \mathcal{I}} B_i) \qquad 
\forall \{B_i: i \in \mathcal{I} \} \subset L^Y, \mbox{ in case } * = 
\wedge$$ 
\end{enumerate}
\end{comments}

\subsection{Injectivity, surjectivity and bijectivity of fuzzy functions} 

\begin{definition}\label{inj}
A fuzzy function $F: (X, E_X) \tto (Y, E_Y)$ is called 
{\it injective}, if 
\begin{itemize}
\item[(inj)]
$F(x,y) * F(x',y) \leq E_X(x,x') \quad 
\forall x,x' \in X, \forall y \in Y$.
\end{itemize}
\end{definition}

\begin{definition}\label{sur}
Given a fuzzy function $F: (X, E_X) \tto (Y, E_Y)$, we define its 
degree of surjectivity by the equality:
$$\sigma(F) := \inf_y \sup_x F(x,y)$$ 
In particular, a fuzzy function $F$ is called $\alpha$-surjective, if 
$\sigma(F) \geq \alpha$.

In case $F$ is injective and $\alpha$-surjective, it is called 
{\it $\alpha$-bijective}.
\end{definition}

\begin{remark}
Let $(X,E_X)$, $(Y,E_Y)$ be $L$-valued sets and $(X',E_{X'})$, 
$(Y',E_{Y'})$ be their subspaces. 
Obviously, the restriction $F'\colon (X',E_{X'}) \tto (Y',E_{Y'})$ 
of an injection $F\colon (X,E_X) \tto (Y,E_Y)$ is an injection. 
The restriction $F'\colon (X,E_{X}) \tto (Y',E_{Y'})$ of an 
$\alpha$-surjection $F\colon (X,E_X) \tto (Y,E_Y)$ is an
$\alpha$-surjection. 
On the other hand, generally the restriction 
$F'\colon (X',E_{X'}) \tto (Y',E_{Y'})$ of an $\alpha$-surjection
$F\colon (X,E_X) \tto (Y,E_Y)$ may fail to be an $\alpha$-surjection.
\end{remark}

A fuzzy function $F: (X, E_X) \tto (Y, E_Y)$ determines a fuzzy 
{\it relation} $F^{-1}: X\times Y \to L$ by setting 
$F^{-1}(y,x) = F(x,y)$\ $\forall x \in X$, $\forall y \in Y$.

\begin{proposition}[Basic properties of injections, $\alpha$-surjections
and $\alpha$-bi\-jections] 
\label{in-sur}
\mbox{}
\begin{enumerate} 
\item 
$F^{-1} :(Y,E_Y) \tto (X,E_X)$ is a fuzzy function iff $F$ is injective
(actually $F^{-1}$ satisfies (3ff) iff F satisfies (inj)) 
\item 
$F$ is $\alpha$-bijective iff $F^{-1}$ is $\alpha$-bijective.
\item 
If $F$ is injective, and L-sets $A_i$ are extensional, then
$$
(\bigwedge_i F^{\rightarrow}(A_i))*(\sigma(F))^2 
\leq 
F^{\rightarrow}(\bigwedge_i A_i) \leq \bigwedge_i F^{\rightarrow}(A_i) 
\qquad \forall \{A_i : i \in {\mathcal{I}}\} \subset L^X.
$$
In particular, if $F$ is $\top$-bijective, then 
$$
(\bigwedge_i F^{\rightarrow}(A_i) = F^{\rightarrow}(\bigwedge_i A_i) 
\qquad \forall \{A_i : i \in {\mathcal{I}}\} \subset L^X.$$
\item 
$F^{\rightarrow}(F^{\leftarrow}(B)) \geq \sigma(F)^2 * B$. 
In particular, if $F$ is $\top$-surjective and $B$ is extensional, then 
$F^{\rightarrow}(F^{\leftarrow}(B)) = B$. 
\item 
$F^{\rightarrow}(c_X) \geq \sigma(F) * c$ where $c_X: X \to L$ is the
constant function with value $c$.
In particular, 
$F^{\rightarrow}(c_X) = c$ if $\sigma(F) = \top$. 
\end{enumerate}
\end{proposition}

\begin{proof}
Properties 1 and 2 follow directly from the definitions.

(2).
The proof of Property 3 is analogous to the proof of item 3 of 
Proposition \ref{im-pr}:

First, reasoning as in the proof of $(\ref{basic1})$ we establish
the following inequality 
\begin{equation}\label{basic3}
\bigvee_{x\in X}\bigl(F(x,y)\bigr)^2 \geq 
\bigl(\bigvee_{x\in X} F(x,y)\bigr)^2.
\end{equation}
In particular, from here it follows that 
\begin{equation}\label{basic4}
\forall y \in Y\ \bigvee_{x\in X} \bigl(F(x,y)^2\bigr) \geq \sigma(F)^2.
\end{equation}
Now, applying $(\ref{basic4})$ and taking into account extensionality of 
$L$-sets $A_i$, we proceed as follows:
$$
\begin{array}{lll}
\multicolumn{3}{l}{
\Bigl(\bigwedge_i F^{\to}(A_i)\Bigr)(x) * \bigl(\sigma(F)\bigr)^2
}
\\
&
\leq& 
\Bigl(\bigwedge_i \bigl(\bigvee_{x_i \in X} (F(x_i,y)*A_i(x_i)\bigr)\Bigr) 
* \bigvee_{x\in X}\bigl(F(x,y)\bigr)^2
\\
&
=&
\bigvee_{x\in X} \Bigl(\bigl(F(x,y)\bigr)^2 * 
\bigwedge_i\bigl(\bigvee_{x_i \in X} F(x_i,y) * A_i(x_i)\bigr)
\\
&
=&
\bigvee_{x\in X}\Bigl(F(x,y) * \bigl(\bigwedge_i \bigl(F(x,y) * 
\bigl(\bigvee_{x_i \in X} F(x_i,y) * A_i(x_i)\bigr)\bigr)\Bigr)
\\
&
=
&
\bigvee_{x\in X}\Bigl(F(x,y) * \bigl(\bigwedge_i \bigl(\bigvee_{x_i\in 
X} (F(x,y) * F(x_i,y)\bigr) \bigr) * A_i(x_i)\bigr)\Bigr)
\\
&
\leq&
\bigvee_{x\in X}\Bigl(F(x,y) * \bigl(\bigwedge_i \bigl(\bigvee_{x_i \in X} 
E_X(x,x_i) * A_i(x_i)\bigr)\bigr)\Bigr)
\\
&
\leq&
\bigvee_{x\in X}\bigl(F(x,y) * \bigl(\bigwedge_i A_i(x)\bigr)\bigr)
\\
&
=&
F^{\to}\bigl(\bigwedge_i A_i)\bigr)(y),
\end{array}
$$
and hence 
$$\Bigl(\bigwedge_i F^{\to}(A_i)\Bigr)* \bigl(\sigma(F)\bigr)^2 
\leq F^{\to}\bigl(\bigwedge_i A_i)\bigr).$$ 
To complete the proof notice that the inequality
$$F^{\rightarrow}(\bigwedge_{i \in \mathcal{I}} A_i) 
\leq 
\bigwedge_{i \in \mathcal{I}} F^{\rightarrow}(A_i)$$ 
is obvious. 

(4).
Let $B \in L^Y$, then 
$$
\begin{array}{lll}
F^{\rightarrow}(F^{\leftarrow}(B))(y)&
=& 
\bigvee_x \bigl(F(x,y)*F^{\gets}(B)(x)\bigr)\\
& 
=&
\bigvee_x 
F(x,y) * \bigl(\bigvee_{y'} F(x,y') * B(y) \bigr)\\
&
\geq&
\bigvee_x F(x,y) * F(x,y) * B(y)\\
&
\geq&
\sigma (F)^2 * B(y),
\end{array}
$$ 
and hence the first inequality in item 4 is proved. From here and 
Proposition \ref{im-pr} (6) the second statement of item 4 follows.

(5)
The proof of the last property is straightforward and therefore omitted.
\end{proof}

\begin{question}
We do not know whether inequalities in items 3, 4 and 5 can be improved.
\end{question}

\begin{comments}
\mbox{}
\begin{enumerate}
\item 
Properties 1 and 2 were first established in \cite{HPS2}.
\item 
In \cite{HPS2} the following version of Property 3 was proved:\\
If $L$ is completely distributive and $F$ is injective, 
then
$$
(\bigwedge_i F^{\rightarrow}(A_i))^5 
\leq 
F^{\rightarrow}(\bigwedge_i A_i) 
\leq 
\bigwedge_i F^{\rightarrow}(A_i) 
\qquad \forall \{A_i : i \in {\mathcal{I}}\} \subset L^X
$$
and
$$F^{\rightarrow}(\bigwedge_i A_i) = 
\bigwedge_i F^{\rightarrow}(A_i) 
\quad \mbox{ in case } \wedge = *$$ 
No extensionality is assumed in these cases.
\item 
In case of an ordinary function $f:X \to Y$ the equality 
$$f^{\to}\bigl(\bigwedge_{i\in \I} (A_i)\bigr) = 
\bigwedge_{i\in\I} f^{\to}(A_i)$$
holds just under assumption that $f$ is injective. 
On the other hand, in case of a fuzzy function $F$ to get a reasonable
counterpart of this property we need to assume that $F$ is bijective. 
The reason for this, as we see it, is that in case of an ordinary
function $f$, when proving the equality, we actually deal only with
points belonging to the image $f(X)$, while the rest of $Y \setminus
f(X)$ does not play any role. 
On the other hand, in case of a fuzzy function $F: X \tto Y$ the whole
$Y$ is ``an image of $X$ to a certain extent'', and therefore, when 
operating with images of $L$-sets, we need to take into account, to what
extent a point $y$ is in the ``image'' of $X$. 
\end{enumerate}
\end{comments}

\section{Further properties of fuzzy category \LFST}

In this section we continue to study properties of the fuzzy category 
\LFST. 
As different from the previous section, were our principal interest 
was in the ``set-theoretic'' aspect of fuzzy functions, here shall 
be mainly interested in their properties of ``categorical nature''.

First we shall specify the two (crisp) categories related to \LFST: 
namely, its bottom frame \LFST$^\bot$ (=\LFST\ (this category was
introduced already in Section 2) and its top frame \LFST$^\top$. 
The last one will be of special importance for us. 
By definition its morphisms $F$ satisfy condition $\mu(F) = \top$, and as
we have seen in the previous section, fuzzy functions satisfying this
condition ``behave themselves much more like ordinary functions'' than
general fuzzy functions. 
Respectively, the results which we are able to establish about
\LFST$^\top$ and about topological category \LFTOP based on it, are more
complete and nice, then their more general counterparts.

Second, note that the ``classical'' category $SET(L)$ of $L$-valued sets
can be naturally viewed as a subcategory of \LFST$^\top$. 
In case $L=\{0,1\}$ obviously the two categories collapse into the
category SET of sets.
On the other hand, starting with the category SET
(=SET$(\{0,1\})$) (i.e.\ $L=\{0,1\}$) of sets and enriching it with
respective fuzzy functions, we obtain again the category SET as
\LFST$^\top$ and obtain the category of sets and partial functions as
\LFST$^\bot$.

\subsection{Preimages of $L$-valued equalities under fuzzy functions}

Let an $L$-valued set $(Y,E_Y)$, a set $X$ and a mapping 
$F: X \times Y \to L$ be given. 
We are interested to find the largest $L$-valued equality $E_X$ on $X$
for which $F: (X,E_X) \to (Y,E_Y)$ is a fuzzy function. 
This $L$-valued equality will be called {\it the preimage of $E_Y$ under
$F$} and will be denoted $F^{\gets}(E_Y)$.

Note first that the axioms
\begin{enumerate}
\item[(1ff)] 
$F(x,y)*E_Y(y,y') \leq F(x,y')$, and
\item[(3ff)] 
$F(x,y)*F(x,y') \leq E(y,y')$ 
\end{enumerate}
do not depend on the $L$-valued equality on $X$ and hence we have to
demand that the mapping $F$ originally satisfies them.
To satisfy the last axiom
\begin{enumerate}
\item[(2ff)] 
$E_X(x,x') * F(x,y) \leq F(x',y)$
\end{enumerate}
in an ``optimal way'' we define 
$$E_X(x,x') := 
\bigwedge_y 
\Bigl(\bigl( F(x,y) \impl F(x',y)\bigr) \wedge 
\bigl(F(x',y) \impl F(x,y) \bigr)\Bigr).$$
Then $E_X: X\times X \to L$ is an $L$-valued equality on $X$. 
Indeed, the validity of properties $E_X(x,x) = \top$ and 
$E_X(x,x') = E_X(x',x)$ is obvious. 
To establish the last property, i.e.\ 
$E_X(x,x') * E_X(x',x'') \leq E_X(x,x'')$, 
we proceed as follows:
$$
\begin{array}{lll}
\multicolumn{3}{l}{
E_X(x,x') * E_X(x',x'')
}
\\
&
=&
\bigwedge_y 
\Bigl(
\bigl(F(x,y) \impl F(x',y)\bigr) \wedge 
\bigl(F(x',y) \impl F(x,y)\bigr)
\Bigr)
\\
&
&
* 
\bigwedge_y 
\Bigl(
\bigl(F(x',y) \impl F(x'',y)\bigr) \wedge 
\bigl(F(x'',y) \impl F(x',y)\bigr)
\Bigr)
\\
&
\leq&
\bigwedge_y 
\Bigl(
\bigl(F(x,y) \impl F(x',y)\bigr) \wedge 
\bigl(F(x',y) \impl F(x,y) \bigr) 
\\
&
&
\ * 
\bigl(F(x',y) \impl F(x'',y)\bigr) \wedge 
\bigl(F(x'',y) \impl F(x',y)\bigr)
\Bigr)
\\
&
\leq&
\bigwedge_y 
\Bigl(
\bigl(F(x,y) \impl F(x',y)\bigr) * 
\bigl(F(x',y) \impl F(x,y) \bigr) 
\\
&
&
\ \wedge
(\bigl(F(x',y) \impl F(x'',y)\bigr) * 
\bigl(F(x'',y) \impl F(x',y) \bigr)
\Bigr)
\\
&
\leq&
\bigwedge_y 
\Bigl(
\bigl(F(x,y) \impl F(x'',y)\bigr) \wedge 
\bigl(F(x'',y) \impl F(x,y)\bigr)
\Bigr)
\\
&
=&
E_X(x,x'').
\end{array}
$$
Further, just from the definition of $E_X$ it is clear that $F$ 
satisfies the axiom {\it (2ff)} and hence it is indeed a fuzzy 
function $F: (X,E_X) \tto (Y,E_Y).$ Moreover, from the definition 
of $E_X$ it is easy to note that it is really the largest $L$-valued 
equality on $X$ for which $F$ satisfies axiom {\it (2ff)}.

Finally, note that the value $\mu(F)$ is an inner property of the 
mapping $F: X\times Y \to L$ and does not depend on $L$-valued 
equalities on these sets.

\begin{question}
We do not know whether the preimage $F^{\gets}(E_Y)$ is the initial 
structure for the source $F: X \tto (Y,E_Y)$ in \LFST. 
Namely, given an $L$-valued set $(Z,E_Z)$ and a ``fuzzy quasi-function''
$G: (Z,E_Z) \tto X$ is it true that composition 
$F \circ G: (Z,E_Z) \tto (Y,E_Y)$ is a fuzzy function if and only if
$G: (Z,E_Z) \tto (X,E_X)$ is a fuzzy function? 
By a fuzzy quasi-function we mean that $G$ satisfies properties {\it
(1ff)} and {\it (3ff)} which do not depend on the equality on $X$.
\end{question}

\subsection{Images of $L$-valued equalities under fuzzy functions}

Let an $L$-valued set $(X,E_X)$, a set $Y$ and a mapping 
$F: X \times Y \to L$ be given. 
We are interested to find the smallest $L$-valued equality $E_Y$ on $Y$
for which $F: (X,E_X) \to (Y,E_Y)$ is a fuzzy function. 
This $L$-valued equality will be called {\it the image of $E_X$ 
under $F$} and will be denoted $F^{\to}(E_X)$.

Note first that the axiom
\begin{enumerate}
\item[(2ff)] 
$E_X(x,x') * F(x,y) \leq F(x',y)$
\end{enumerate}
does not depend on the $L$-valued equality on $Y$ and hence we have to
demand that the mapping $F$ originally satisfies it.
Therefore we have to bother that $F$ satisfies the remaining two axioms:
\begin{enumerate}
\item[(1ff)] $F(x,y)*E_Y(y,y') \leq F(x,y')$, and
\item[(3ff)] $F(x,y)*F(x,y') \leq E(y,y')$ 
\end{enumerate}
These conditions can be rewritten in the form of the double inequality:
$$
\begin{array}{lll}
F(x,y)*F(x,y')&
\leq&
E_Y(y,y')\\
&
\leq&
\bigl(F(x,y')\impl F(x,y)\bigr) \wedge \bigl(F(x,y) \impl F(x,y')\bigr).
\end{array}
$$
Defining $E_Y$ by the equality
$$E_Y(y,y') = \bigvee_x \bigl(F(x,y) * F(x,y')\bigr),$$
we shall obviously satisfy both of them.
Moreover, it is clear that $E_Y$ satisfies property (3ff) and besides
$E_Y$ cannot be diminished without loosing this property. 
Hence we have to show only that $E_Y$ is indeed an $L$-valued 
equality. 
However, to prove this we need the assumption that $\sigma(F) = \top$,
that is $F$ is $\top$-surjective.
Note that 
$$E_Y(y,y) = \bigvee_x \bigl(F(x,y)*F(x,y)\bigr) \geq (\sigma(F))^2,$$ 
and hence the first axiom is justified in case $\sigma(F) = \top$. 

The equality $E_Y(y,y') = E_Y(y',y)$ is obvious.

Finally, to establish the last property, we proceed as follows. 
Let $y,y',y'' \in Y$. 
Then
$$
\begin{array}{lll}
\multicolumn{3}{l}{
E_Y(y,y')*E_Y(y',y'')
}
\\
&
=&
\bigvee_x\bigl(F(x,y)*F(x,y')\bigr)
\,
*
\, 
\bigvee_x\bigl(F(x,y')*F(x,y'')\bigr)
\\
&
=&
\bigvee_{x,x'} \bigl(F(x,y)*F(x,y')* 
F(x',y')*F(x',y'')\bigr)\\
&
\leq&
\bigvee_{x,x'} \bigl(F(x,y)* E_X(x,x')*F(x',y'')\bigr)\\
&
\leq&
\bigvee_x \bigl(F(x,y) * F(x,y'')\bigr)
\\
&
=&
E(y,y'').
\end{array}
$$

\begin{question}
We do not know whether the image $F^{\to}(E_X)$ is the final 
structure for the sink $F: (X,E_X) \tto Y$ in \LFST\ in case 
$\sigma(F) = \top$. 
Namely, given an $L$-valued set $(Z,E_Z)$ and a ``fuzzy 
almost-function'' $G: Y \tto (Z,E_Z) $ is it true that composition 
$F \circ G: (Z,E_Z) \tto (Y,E_Y)$ is a fuzzy function if and only if
$G: (Y,E_Y) \tto (Z,E_Z)$ is a fuzzy function? 
By a fuzzy almost-function we mean that $G$ satisfies property
(2ff) which does not depend on the equality on $Y$.
\end{question}

\subsection{Products in \LFSTT}

Let $\Y = \{(Y_i, E_i): i \in \I \}$ be a family of $L$-valued sets and
let $Y = \prod_i Y_i $ be the product of the corresponding sets. 
We introduce the $L$-valued equality $E: Y\times Y \to L$ on $Y$ by
setting $E_Y(y,y') = \bigwedge_{i\in\I} E_i(y_i,{y'}_i)$ where
$y=(y_i)_{i\in\I}$, $y'=({y'}_i)_{i\in\I}$.
Further, let $p_i :Y \to Y_i$ be the projection. 
Then the pair $(Y,E)$ thus defined with the family of projections 
$p_i :Y \to Y_i$, $i \in \I$, is the product of the family $\Y$ in the
category \LFSTT.

To show this notice first that, since the morphisms in this category are 
fuzzy functions, a projection $p_{i_0}: Y \to Y_{i_0}$ must be realized 
as the fuzzy function $p_{i_0}: Y \times Y_{i_0} \to L$ such that 
$p_{i_0}(y,y^0_{i_0}) = \top$ if and only if the $i_0$-coordinate 
of $y$ is $y^0_{i_0}$ and $p_{i_0}(y,y^0_{i_0}) = \bot$ otherwise.

Next, let $F_i: (X,E_X) \tto (Y_i,{E_i})$, $i \in \I$ be a family of
fuzzy functions. 
We define the fuzzy function $F: (X,E_X) \tto (Y,E_Y)$ by the equality:
$$F(x,y) = \bigwedge_{i\in\I} F_i(x,y_i).$$
It is obvious that $\mu(F) = \top$ and hence $F$ is in \LFSTT.
Finally, notice that the composition
$$(X,E_X) \ 
\stackrel{F}{\longrightarrow} \ 
(Y,E_Y) \
\stackrel{p_{i_0}}{\longrightarrow} \ 
(Y_{i_0},E_{i_0})$$
is the fuzzy function 
$$F_{i_0}: (X,E_X) \tto (Y_{i_0},E_{i_0}).$$
Indeed, let $x^0 \in X$ and $y^0_{i_0} \in Y_{i_0}$. 
Then, taking into account that $\mu(F_i) = \top$ for all $i \in \I$, 
we get
$$
\begin{array}{lll}
(p_{i_0} \circ F)(x^o,y^0_{i_0})&
=&
\bigvee_{y\in Y} 
\bigl(p_{i_0}(y,y^0_{i_0}) \wedge F(x^0,y)\bigr)\\
&
=&
\bigvee_{y\in Y} 
\bigl(p_{i_0}(y,y^0_{i_0}) \wedge \bigwedge_{i\in\I} F_i(x^0,y_i)\bigr)\\
&
=&
F_{i_0}(x^0,y_{i_0}).
\end{array}$$

\begin{question}
We do not know whether products in \LFST\ can be defined in a reasonable
way.
\end{question}

\subsection{Coproducts in \LFST}

Let $\X$ = $\{(X_i,E_i): i \in \I \}$ be a family of $L$-valued sets,
let $X = \bigcup X_i$ be the disjoint sum of sets $X_i$.
Further, let $q_i: X_i \to X$ be the inclusion map. 
We introduce the $L$-equality on $X_0$ by setting $E(x,x') = E_i(x,x')$
if $(x,x') \in X_i \times X_i$ for some $i \in \I$ and $E(x,x') = \bot$
otherwise (cf.\ \cite{Ho92}). 
An easy verification shows that $(X,E)$ is the coproduct of $\X$ in 
\LFST\ and hence, in particular, in \LFSTT. 

Indeed, given a family of fuzzy functions $F_i: (X_i,E_i) \to (Y,E_Y)$,
let the fuzzy function 
$$\oplus_{i\in\I} F_i: (X,E) \to (Y,E_Y)$$
be defined by
$$\oplus_{i\in\I} F_i(x,y) = F_{i_0}(x,y) \mbox{ whenever } x \in X_{i_0}.$$
Then for $x=x_{i_0} \in X_{i_0}$ we have
$$
\begin{array}{lll}
\bigl(\oplus_{i\in\I} F_i \circ q_{i_0} \bigr)(x,y)&
=&
\bigvee_{x'\in X}
\Bigl(q_{i_0}(x,x') \wedge \bigl(\oplus_{i\in\I} F_i (x',y)\bigr)\Bigr)\\
&
=&
F_{i_0}(x,y).
\end{array}
$$

\subsection{Subobjects in \LFST}

Let $(X,E)$ be an $L$-valued set, let $Y \subset X$ and let $e: Y \to X$
be the natural embedding. 
Further, let $E_Y := e^{\gets}(E)$ be the preimage of the $L$-valued
equality $E$. 
Explicitly, in this case this means that $E_Y(y,y') = E(y,y')$ for all
$y,y' \in Y$.
One can easily see that $(Y,E_Y)$ is a subobject of $(X,E)$ in the fuzzy 
category \LFST.

\section{Fuzzy category \LFTop} 
 
\subsection{Basic concepts}

\begin{definition}[see \cite{So2000}, cf.\ also \cite{Ch}, \cite{Go73}, 
\cite{HoSo99}]\label{LFTop} 
A family $\tau_X \subset L^X$ of {\it extensional}
$L$-sets\footnote{Since $L$-topology is defined on an {\it $L$-valued set} 
$X$ the condition of extensionality of elements of $L$-topology seems 
natural. 
Besides the assumption of extensionality is already implicitly included
in the definition of a fuzzy function.} 
is called an $L$-topology on an $L$-valued set $(X,E_X)$ if it is closed
under finite meets, arbitrary joins and contains $0_X$ and $1_X$. 
Corresponding triple $(X,E_X,\tau_X)$ will be called an $L$-valued
$L$-topological space or just an $L$-topological space for short. 
A fuzzy function $F: (X,E_X,\tau_X) \tto (Y,E_Y,\tau_Y)$ is called 
{\it continuous} if $F^{\leftarrow}(V) \in \tau_X$ for all 
$V \in \tau_Y$. 
\end{definition} 

$L$-topological spaces and continuous fuzzy mappings between them form
the fuzzy category which will be denoted \LFTop.
Indeed, let 
$$F\colon (X,E_X,\tau_X) \tto (Y,E_Y,\tau_Y)\
\mbox{and}\
G\colon (Y,E_Y,\tau_Y) \tto (Z,E_Z,\tau_Z)$$
be continuous fuzzy functions and let $W \in \tau_Z$. 
Then 
$$
\begin{array}{lll}
(G \circ F)^{\gets}(W)(x)&
=&
\bigvee_z\Bigl(\bigl(G \circ F)(x,z) * W(z)\Bigr)\\
&
=&
\bigvee_{z,y} \Bigl(F(x,y)*G(y,z)*W(z)\Bigr).
\end{array}
$$
On the other hand, $G^{\gets}(W)(y) = \bigvee_z G(y,z) * W(z)$ and 
$$
\begin{array}{lll}
F^{\gets}\bigl(G^{\gets}(W)\bigr)(x)&
=&
\bigvee_y F(x,y) *\bigl(\bigvee_z (G(y,z) * W(z))\bigr)\\
&
=& 
\bigvee_{z,y} \Bigl(F(x,y)*G(y,z)*W(z)\Bigr).
\end{array}
$$
Thus $(G \circ F)^{\gets}(W) = G^{\gets}(F^{\gets}(W))$ for every $W$,
and hence composition of continuous fuzzy functions is continuous.
Besides, we have seen already before that 
$\mu(G \circ F) \geq \mu(G) * \mu(F)$.
Finally, $E_X^{\leftarrow}(B) = B$ for every {\it extensional} 
$B \in L^X$ and hence the identity mapping 
$E_X: (X,E_X,\tau_X) \tto (X,E_X,\tau_X)$ is continuous.

\begin{remark}
In case when $L$-valued equality $E_X$ is crisp, i.e.\ when $X$ is an 
ordinary set, the above definition of an L-topology on $X$ reduces to the 
``classical'' definition of an $L$-topology in the sense of Chang and 
Goguen, \cite{Ch}, \cite{Go73}. 
\end{remark}

\begin{remark}
Some (ordinary) subcategories of the fuzzy category \LFTop\ will be of 
special interest for us. 
Namely, let \LFTop$^\bot$ =: FTOP(L) denote the bottom frame of \LFTop,
let \LFTOP\ be the top frame of \LFTop, and finally let L-TOP(L) denote
the subcategory of \LFTop\ whose morphisms are ordinary functions. 
Obviously the ``classical'' category L-TOP of Chang-Goguen $L$-topological
spaces can be obtained as a full subcategory L-TOP(L) whose objects carry 
crisp equalities. 
Another way to obtain L-TOP is to consider fuzzy subcategory of \LFTop
whose objects carry crisp equalities and whose morphisms satisfy
condition $\mu(F) > \bot.$
\end{remark}

In case when $L$ is an $MV$-algebra and involution $^c: L \to L$ on $L$
is defined in the standard way, i.e.\ $\alpha^c := \alpha \impl \bot$ 
we can reasonably introduce the notion of a closed $L$-set in an
$L$-topological space:

\begin{definition}
An $L$-set $A$ in an $L$-topological space $(X,E_X,\tau_X)$ is called 
closed if $A^c \in \tau_X$ where $A^c \in L^X$ is defined by the equality
$$A^c(x) := A(x) \impl \bot \quad \forall x \in X.$$
\end{definition}

Let $\C_X$ denote the family of all closed $L$-sets in $(X,E_X,\tau_X)$.
In case when $L$ is an $MV$-algebra the families of sets $\tau_X$ and 
$\C_X$ mutually determine each other:
$$A \in \tau_X \Longleftrightarrow A^c \in \C_X.$$

\subsection{Analysis of continuity}

Since the operation of taking preimages $F^{\leftarrow}$ commutes with 
joins, and in case when $\mu(F)=\top$ also with meets (see Proposition 
\ref{im-pr}), one can easily verify the following

\begin{theorem}\label{cont}
Let $(X,E_X,\tau_X)$ and $(Y,E_Y,\tau_Y)$ be $L$-topological spaces, 
$\beta_Y$ be a base of $\tau_Y$, $\xi_Y$ its subbase and let $F: X \tto Y$ 
be a fuzzy function. 
Then the following are equivalent:
\begin{enumerate}
\item[(1con)] 
$F$ is continuous;
\item[(2con)] 
for every $V \in \beta_Y$ it holds $F^{\leftarrow}(V) \in \tau_X$;
\item[(3con)] 
$F^{\leftarrow}(Int_Y(B)) \leq Int_X (F^{\leftarrow}(B))$, 
for every $B \in L^Y$ where $Int_X$ and $Int_Y$ are the corresponding
$L$-interior operators on $X$ and $Y$ respectively. 
\end{enumerate}
In case when $\mu(F) = \top$ these conditions are equivalent also to the 
following
\begin{enumerate}
\item[(4con)] 
for every $V \in \xi_Y$ it holds $F^{\leftarrow}(V) \in \tau_X$.
\end{enumerate}
\end{theorem}

In case when $L$ is an MV-algebra one can characterize continuity of a 
fuzzy function by means of closed $L$-sets and $L$-closure operators:

\begin{theorem}\label{cont-cl}
Let $(L,\leq,\vee,\wedge,*)$ be an MV-algebra, $(X,E_X,\tau_X)$ and 
$(Y,E_Y,\tau_Y)$ be $L$-topological spaces and $F: X \tto Y$ be a 
fuzzy function. 
Further, let $\C_X$, $\C_Y$ denote the families of closed $L$-sets and
$cl_X$, $cl_Y$ denote the closure operators in $(X,E_X,\tau_X)$ and 
$(Y,E_Y,\tau_Y)$ respectively. 
Then the following two conditions are equivalent:
\begin{enumerate}
\item[(1con)] 
$F$ is continuous;
\item[(5con)] 
For every $B \in \C_Y$ it follows $F^{\leftarrow}(B) \in \C_X$.
\end{enumerate}
In case when $\mu(F) = \top$, the previous conditions are equivalent to
the following:
\begin{enumerate}
\item[(6con)] 
For every $A \in L^X$ it holds $F^{\rightarrow}(cl_X(A)) \leq 
cl_Y(F^{\rightarrow}(A)).$
\end{enumerate}
\end{theorem}

\begin{proof}
In case when $L$ is equipped with an order reversing involution, as it is
in our situation, families of closed and open $L$-sets mutually determine
each other. 
Therefore, to verify the equivalence of (1con) and (5con) it is sufficient
to notice that for every $B \in L^Y$ and every $x \in X$ it holds
$$
\begin{array}{lll}
F^{\leftarrow}(B^c)(x)&
=&
\bigvee_y \bigl(F(x,y) * (B(y) \impl \bot)\bigr)\\
&
=&
\bigvee_y \bigl(F(x,y) * B(y) \impl \bot \bigr)\\
&
=&
\bigl(\bigvee_y (F(x,y) * B(y))\bigr)^c\\
&
=&
(F^{\leftarrow}(B))^c(x),
\end{array}
$$
and hence
$$
F^{\leftarrow}(B^c) = ( F^{\leftarrow}(B))^c \quad \forall B \in L^Y,
$$
i.e.\ operation of taking preimages preserves involution.

To show implication (5con) $\Longrightarrow$ (6con) under assumption 
$\mu(F)=\top$ let $A \in L^X$. 
Then, according to Proposition \ref{im-pr} (5), 
$$A \leq F^{\leftarrow}(F^{\rightarrow}(A)) 
\leq F^{\leftarrow}(cl_Y(F^{\rightarrow}(A))),$$
and hence, by (5con), also 
$$cl_X (A) \leq F^{\leftarrow}(cl_Y (F^{\rightarrow}(A))).$$ 
Now, by monotonicity of the image operator and by Proposition \ref{im-pr}
(6) (taking into account that $cl_X A$ is extensional as a closed
$L$-set), we get: 
$$ 
F^{\rightarrow}(cl_X (A)) 
\leq 
F^{\rightarrow}\bigl(F^{\leftarrow}(cl_Y (F^{\rightarrow}(A)))\bigr) 
\leq 
cl_Y (F^{\rightarrow}(A)).$$

Conversely, to show implication (6con) $\Longrightarrow$ (5con) 
let $B \in \C_Y$ and let $F^{\gets}(B) := A$. Then, by (6con), 
$$F^{\to}(cl_X (A)) \leq cl_Y (F^{\to} (A)) \leq cl_Y (B) = B.$$
In virtue of Proposition \ref{im-pr} (5) and taking into account that
$\mu(F) = \top$, it follows from here that 
$cl_X (A) \leq F^{\gets}(B) = A$, and hence $cl_X (A) = A$.
\end{proof}

\subsection{Fuzzy $\alpha$-homeomorphisms and fuzzy $\alpha$-homeomorphic 
spaces}

The following definition naturally stems from Definitions \ref{inj}, 
\ref{sur} and \ref{LFTop} and item 2 of Proposition \ref{im-pr}:

\begin{definition} 
Given $L$-topological spaces $(X,E_X,\tau_X)$ and $(Y,E_Y,\tau_Y)$, a
fuzzy function $F: X \tto Y$ is called a fuzzy $\alpha$-homeomorphism if
$\mu(F) \geq \alpha$, $\sigma(F) \geq \alpha$, it is injective, 
continuous, and the inverse fuzzy function $F^{-1}: Y \tto X$ is also
continuous.
Spaces $(X,E_X,\tau_X)$ and $(Y,E_Y,\tau_Y)$ are called fuzzy 
$\alpha$-homeomorphic if there exists a fuzzy $\alpha$-homeomorphism
$F: (X,E_X,\tau_X) \tto (Y,E_Y,\tau_Y)$.
\end{definition}

One can easily verify that composition of two fuzzy 
$\alpha$-homeomorphisms is a fuzzy $\alpha^2$-home\-omorphism; 
in particular, composition of fuzzy $\top$-homeo\-morphisms is a fuzzy 
$\top$-homeo\-morphism, and hence fuzzy $\top$-homeo\-morph\-isms
determine the equivalence relation $\stackrel{\top}{\approx}$ on the class
of all $L$-topological spaces. 
Besides, since every (usual) homeomorphism is 
obviously a fuzzy $\top$-homeo\-morphism, homeomorphic spaces are also 
fuzzy $\top$-homeo\-morphic: 
$$(X,E_X,\tau_X) \approx (Y,E_Y,\tau_Y) 
\Longrightarrow (X,E_X,\tau_X) \stackrel{\top}{\approx} (Y,E_Y,\tau_Y).$$ 
The converse generally does not hold:

\begin{example}
Let $L$ be the unit interval $[0,1]$ viewed as an $MV$-algebra 
(i.e.\ $\alpha * \beta = max\{\alpha+\beta-1, 0\}$, let $(X,\varrho)$ be
an uncountable separable metric space such that $\varrho(x,x') \leq 1$
$\forall x,x' \in X$, and let $Y$ be its countable dense subset. 
Further, let the $L$-valued equality on $X$ be defined by 
$E_X(x,x') := 1 - \varrho(x,x')$ and let $E_Y$ be its restriction to $Y$. 
Let $\tau_X$ be any $L$-topology on an $L$-valued set $(X,E_X)$ (in
particular, one can take $\tau_X := \{c_X \mid c \in [0,1] \}$). 
Finally, let a fuzzy function $F: (X,E_X) \tto (Y,E_Y)$ be defined by
$F(x,y) := 1 - \varrho(x,y)$. 
It is easy to see that $F$ is a $\top$-homeomorphism and hence 
$(X,E_X,\tau_X) \stackrel{\top}{\approx} (Y,E_Y,\tau_Y)$. 
On the other hand $(X,E_X,\tau_X) \not{\approx} (Y,E_Y,\tau_Y)$ just for
set-theoretical reasons.
\end{example}

\section{Category \LFTOP}

Let \LFTOP\ be the top-frame of \LFTop, i.e.\ \LFTOP\ is a category
whose objects are the same as in \LFTop, that is L-topological spaces,
and morphisms are continuous fuzzy functions 
$F: (X,E_X,\tau_X) \tto (Y,E_Y,\tau_Y)$ such that $\mu(F) = \top.$

Note, that as different from the fuzzy category \LFTop, \LFTOP\ is a 
usual category. Applying Theorem \ref{cont}, we come to the following 
result:

\begin{theorem}\label{cont-top}
Let $(X,E_X,\tau_X)$ and $(Y,E_Y,\tau_Y)$ be $L$-topological spaces, 
$\beta_Y$ be a base of $\tau_Y$, $\xi_Y$ its subbase and let $F: X \tto Y$ 
be a fuzzy function. 
Then the following conditions are equivalent:
\begin{enumerate}
\item[(1con)] 
$F$ is continuous;
\item[(2con)] 
for every $V \in \beta_Y$ it holds $F^{\leftarrow}(V) \in \tau_X$;
\item[(3con)] 
$F^{\leftarrow}(Int_Y(B)) \leq Int_X (F^{\leftarrow}(B))$, where
$Int_X$ and $Int_Y$ are the corresponding $L$-interior operators on $X$
and $Y$ respectively;
\item[(4con)] 
for every $V \in \xi_Y$ it holds $F^{\leftarrow}(V) \in \tau_X$.
\end{enumerate}
\end{theorem}

In case when $L$ is an MV-algebra, we get from \ref{cont-cl}

\begin{theorem}\label{cont-cltop}
Let $(L,\leq,\vee,\wedge,*)$ be an MV-algebra, $(X,E_X,\tau_X)$ and 
$(Y,E_Y,\tau_Y)$ be $L$-topological spaces and $F: X \tto Y$ be a 
a morphism in \LFTOP. 
Further, let $\C_X$, $\C_Y$ denote the families of closed $L$-sets and
$cl_X$, $cl_Y$ denote the closure operators on $(X,E_X,\tau_X)$ and 
$(Y,E_Y,\tau_Y)$ respectively. 
Then the following two conditions are equivalent:
\begin{enumerate}
\item[(1con)] 
$F$ is continuous;
\item[(5con)] 
For every $B \in \C_Y$ it follows $F^{\leftarrow}(B) \in \C_X$;
\item[(6con)] 
For every $A \in L^X$ it holds 
$F^{\rightarrow}(cl_X(A)) \leq cl_Y(F^{\rightarrow}(A))$.
\end{enumerate}
\end{theorem}

\begin{theorem}\label{topcat}
$$\mbox{
\LFTOP is a topological category over the category \LFSTT.
}$$
\end{theorem}

\begin{proof}
Since intersection of any family of L-topologies is an L-topology, 
\LFTOP\ is fiber complete. 
Therefore we have to show only that any structured source in \LFSTT\
$F: (X,E_X) \tto (Y,E_Y,\tau_Y)$ has a unique initial lift. 
Let 
$$\tau_X := F^{\gets}(\tau_Y):= \{F^{\gets}(V) \mid V \in \tau_Y \}.$$
Then from theorem \ref{im-pr} it follows that $\tau_X$ is closed under 
taking finite meets and arbitrary joins. 
Furthermore, obviously $F^{\gets}(0_Y) = 0_X$ and taking into account
condition $\mu(F)=\top$ one easily establishes that $F^{\gets}(1_Y) = 1_X$. 
Therefore, taking into account that preimages of extensional $L$-sets are 
extensional, (see Subsection \ref{impref}) we conclude that the family
$\tau_X$ is
an $L$-topology on $X$. 

Further, just from the construction of $\tau_X$ it is clear that 
$F: (X,E_X,\tau_X) \tto (Y,E_Y,\tau_Y)$ is continuous and, moreover,
$\tau_X$ is the weakest L-topology on $X$ with this property. 

Let now $(Z,E_Z,\tau_Z)$ be an $L$-topological space and 
$H: (Z,E_Z) \tto (X,E_X)$ a fuzzy function such that the composition 
$G := H \circ F: (Z,E_Z,\tau_Z) \tto (Y,E_Y,\tau_Y)$ is continuous. 
To complete the proof we have to show that $H$ is continuous. 

Indeed, let $U \in \tau_X$. 
Then there exists $V \in \tau_Y$ such that $U = F^{\gets}(V)$. 
Therefore 
$$H^{\gets}(U) = H^{\gets}\bigl(F^{\gets}(V)\bigr) 
= G^{\gets}(V) \in \tau_Z$$
and hence $H$ is continuous.
\end{proof}

\subsection{Products in \LFTOP}
 
Our next aim is to give an explicite description of the product in \LFTOP. 

Given a family $\Y = \{(Y_i, E_i, \tau_i): i \in \I \}$ of 
$L$-topological spaces, let $(Y,E)$ be the product of the corresponding
$L$-valued sets $\{(Y_i, E_i): i \in \I \}$ in \LFSTT and let 
$p_{i}: Y \to Y_{i}$ be the projections. 

Further, for each $U_i \in \tau_i$ let $\hat U_i := p_i^{-1}(U_i)$.
Then the family $\xi := \{\hat U_i: U_i \in \tau_i, i \in \I \}$ is a 
subbase of an $L$-topology $\tau$ on the product $L$-valued set $(X,E_X)$
which is known to be the product $L$-topology for L-topological spaces
$\{(Y_i,\tau_i) \mid i \in \I \}$ in the category L-TOP. 
In its turn the triple $(Y,E,\tau)$ is the product of $L$-topological
spaces $\{(Y_i, E_i, \tau_i): i \in \I \}$ in the category \LFTOP. 
Indeed, let $(Z,E_Z,\tau_Z)$ be an $L$-topological space and 
$\{ F_i: Z \tto Y_i \mid i \in \I \}$ be a family of continuous fuzzy
mappings. 
Then, defining a mapping $F: Z\times Y \to L$ by 
$F(z,y) = \wedge_{i\in\I} F_i(z,y_i)$ we obtain a fuzzy function 
$F: Z \tto Y$ such that $\mu(F) = \bigwedge_{i\in\I}\mu(F_i) = \top$ and
besides it can be easily seen that for every $i_0 \in \I$, every $z \in
Z$, and for every $U_{i_0} \in \tau_{i_0}$ it holds:
$$\begin{array}{lll}
F^{\leftarrow}(\hat U_{i_0})(z)&
=& 
\bigvee_{y \in Y}
\Bigl(\bigl(\bigwedge_{i\in\I}F_i(z,y_i)\bigr) * U_{i}(y_i)\Bigr)
\\ 
&
=& 
\bigvee_{y_i \in Y_i}
\Bigl(\bigwedge_{i\in\I} \bigl(F_i(z,y_i)\bigr) * U_{i_0}(y_{i_0})\Bigr)\\
&
=&
\bigwedge_{i\in \I}
\Bigl(\bigvee_{y_i \in Y_i} F_i(z,y_i)\Bigr) 
\wedge 
\bigvee_{y_{i_0} \in X_{i_0}}
\bigl(F_{i_0} (z,y_{i_0}) * U_{i_0}(y_{i_0})\bigr)
\\
&
=&
\top \wedge F^{\gets}_{i_0}(U_{i_0})(z)\\
&
=&
F^{\gets}_{i_0}(U_{i_0})(z).
\end{array}$$ 
Hence continuity of all $F_i$ guarantees the continuity of 
$F: (Z,E_Z,\tau_Z) \tto (Y,E,\tau)$. 
Thus $(Y,E,\tau)$ is indeed the product of 
$\Y = \{(Y_i, E_i, \tau_i): i \in \I \}$ in \LFST$^\top$.

\begin{question}
We do not know whether products in the fuzzy category \LFTop\ exist.
\end{question}

\subsection{Subspaces in \LFTOP}

Let $(X,E,\tau)$ be an $L$-valued $L$-topological space, $Y \subset X$ and 
let $(Y,E_Y)$ be the subobject of the L-valued set $(X,E)$ in the category 
\LFST. 
Further, let $\tau_Y$ be the subspace L-topology, that is $(Y,\tau_Y)$ is
a subspace of the $L$-topological space $(X,\tau)$ in the category
L-TOP. 
Then it is clear that the triple $(Y,E_Y,\tau_Y)$ is the subobject of
$(X,E, \tau)$ in the category \LFTOP (and in the fuzzy category \LFTop\
as well).

\subsection{Coproducts in \LFTOP}

Given a family $\X := \{(X_i,E_i,\tau_i) \}$ of $L$-topological spaces
let $(X,E)$ be the direct sum of the corresponding $L$-valued sets 
$(X_i,E_i)$ in \LFST. 
Further, let $\tau$ be the $L$-topology on $X$ determined by the subbase
$\bigcup_{i\in\I} \tau_i \subset 2^X$. 
In other words $(X,\tau)$ is the coproduct of $L$-topological spaces
$(X_i,\tau_i)$ in the category L-TOP. 
Then the triple $(X,E,\tau)$ is the coproduct of the family 
$\X := \{(X_i,E_i,\tau_i) \}$ in the category \LFTOP (and in the fuzzy
category \LFTop as well). 
Indeed, let $q_i: (X_i,E_i,\tau_i) \to (X,E,\tau), i\in\I$ denote the
canonical embeddings. 
Further, consider an $L$-topological space $(Y,E_Y,\tau_Y)$ and a family
of continuous fuzzy functions $F_i: (X_i,E_i,\tau_i) \tto (Y,E_Y,\tau_Y)$. 
Then, by setting $F(x,y) := F_i(x_i,y)$ whenever $x = x_i \in X_i$, we
obtain a continuous fuzzy function $F: (X,E,\tau) \tto (Y,E_Y,\tau_Y)$
(i.e.\ a mapping $F: X\times Y \to L$) such that $F_i = q_i \circ F$ for
every $i \in \I$.

\subsection{Quotients in \LFTOP}

Let $(X,E_X,\tau_X)$ be an $L$-topological space, 
let $q: X \to Y$ be a surjective mapping. 
Further, let $q^{\to}(E_X) =: E_Y$ be the image of the $L$-valued
equality $E_X$ and let
$\tau_Y = \{V \in L^Y \mid q^{-1}(V) \in \tau_X \}$, that is $\tau_Y$ is
the quotient $L$-topology determined by the mapping $q: (X,\tau) \to Y$
in the category L-TOP. 
Then $(Y,E_Y,\tau_Y)$ is the quotient object in the category \LFTOP. 
Indeed, consider a fuzzy function $F: (X,E_X,\tau_X) \tto (Z,E_Z,\tau_Z)$
and let $G: (Y,E_Y) \tto (Z,E_Z)$ be a morphism in \LFST\ such that 
$q \circ G = F$. 
Then an easy verification shows that the fuzzy function 
$G: (Y,E_Y,\tau_Y) \tto (Z,E_Z,\tau_Z)$ is continuous (i.e.\ a morphism
in \LFTOP) if and only if $F: (X,E_X,\tau_X) \to (Z,E_Z,\tau_Z)$ is 
continuous (i.e.\ a morphism in \LFTOP).

Our next aim is to consider the behaviour of some topological properties
of L-valued L-topological spaces in respect of fuzzy function. In this
work we restrict our interest to the property of compactness. Some other
topological properties, in particular, connectedness and separation
properties will be studied in a subsequent work.

\section{Compactness} 

\subsection{Preservation of compactness by fuzzy functions}

One of the basic facts of general topology --- both classic and 
``fuzzy'', is preservation of compactness type properties by continuous
mappings. 
Here we present a counterpart of this fact in \LFTop. 
However, since in literature on fuzzy topology different definitions of
compactness can be found, first we must specify which one of compactness
notions will be used. 

\begin{definition}\label{comp}
An $L$-topological space $(X,E,\tau)$ will be called 
$(\alpha, \beta)$-compact where $\alpha, \beta \in L$, if for every family 
$\U \subset \tau$ such that $\bigvee \U \geq \alpha$ there exists a finite
subfamily $\U_0 \subset \U$ such that $\bigvee\U_0 \geq \beta$.
An $(\alpha,\alpha)$-compact space will be called just 
$\alpha$-compact.\footnote{Note that Chang's definition of compactness
\cite{Ch} for a $[0,1]$-topological space is equivalent to our 
$1$-compactness.
An $[0,1]$-topological space is compact in Lowen's sense \cite{Lo76} if it 
is $(\alpha,\beta)$-compact for all $\alpha \in [0,1]$ and all 
$\beta < \alpha$}
\end{definition}

\begin{theorem}\label{pr-comp}
Let $(X,E_X, \tau_X)$, $(Y, 
E_Y, \tau_Y)$ be $L$-topological spaces, $F:X \tto Y $ be a continuous fuzzy 
function such that $\mu(F) \geq \beta$, and $\sigma(F) \geq \gamma$. If $X$ 
is $\alpha * \beta$-compact, then $Y$ is $(\alpha, 
\alpha*\beta*\gamma)$-compact.
\end{theorem}

\begin{proof}
Let $\V \subset \tau_Y$ be such that $\bigvee \V \geq \alpha$. 
Then, applying Proposition \ref{im-pr} (4), (7) and taking in view 
monotonicity of $F^{\gets}$, we get 
$$\bigvee_{V \in \V} F^{\leftarrow}(V) = 
F^{\leftarrow}(\bigvee_{V \in \V}V) \geq F^{\leftarrow}(\alpha) \geq \alpha 
* \beta.$$ Now, since $(X,E_X,\tau_X)$ is $\alpha * \beta$-compact, it 
follows that there exists a finite subfamily $\V_0 \subset \V$ 
such that 
$$\bigvee_{V\in 
\V_0} F^{\leftarrow}(V) \geq \alpha*\beta.$$ 
Applying Propositions 
\ref{im-pr} (6),(4) and \ref{in-sur} (5) we obtain:
$$\bigvee_{V\in\V_0} V \geq 
F^{\rightarrow}\Bigl(F^{\leftarrow}\bigl(\bigvee_{V\in\V_0} 
V\bigr)\Bigr) =
F^{\rightarrow}\Bigl(\bigvee_{V\in\V_0} 
\bigl(F^{\leftarrow}(V)\bigr)\Bigr)
\geq 
F^{\rightarrow}\bigl(\alpha * \beta\bigr) \geq 
\alpha*\beta*\gamma.$$
\end{proof}

\begin{corollary} Let $(X,E_X, \tau_X)$, $(Y,E_Y,\tau_Y)$ be 
$L$-topological spaces, $F:X \tto Y$ be a fuzzy function such that $\mu(F) 
= \top $ and $\sigma(F) = \top$. If $X$ is $\alpha$-compact, 
then $Y$ is also $\alpha$-compact.
\end{corollary}

\subsection{Compactness in case of an $MV$-algebra}

In case $L$ is an MV-algebra one can characterize compactness by systems
of closed $L$-sets:

\begin{proposition}
Let $(X,E_X,\tau_X)$ be an $L$-topological space and let $\C_X$ be the 
family of its closed $L$-sets. 
Then the space $(X,E_X,\tau_X)$ is $(\alpha,\beta)$-compact if and only
if for every $\A \subset \C_X$ the following implication follows:
$$
\mbox{if } \bigwedge_{A\in \A_0} A \not\leq \beta^c 
\mbox{ for every finite family } \A_0 \subset \A, 
\mbox{ then } \bigwedge_{A\in \A} A \not\leq \alpha^c.
$$
\end{proposition}

\begin{proof}
One has just to take involutions ``$\impl \bot$'' in the definition of 
$(\alpha,\beta)$-compactenss and apply De Morgan law.
\end{proof}

\subsection{Perfect mappings: case of an MV-algebra $L$}

In order to study preservation of compactness by preimages of fuzzy 
functions we introduce the property of $(\alpha,\beta)$-perfectness of a 
fuzzy function. 
Since we shall operate with closed $L$-sets, from the beginning it will
be assumed that $L$ is an $MV$-algebra. 

First we shall extend the notion of compactness for $L$-subsets of 
$L$-topolog\-ical spaces.
We shall say that an $L$-set $S: X \to L$ is $(\alpha,\beta)$-compact if 
for every family $\A$ of closed $L$-sets of $X$ the following 
implication holds:
$$
\mbox{ if } S\wedge\bigl(\bigwedge_{A\in \A_0} A\bigr) \not\leq \beta^c 
\mbox{ for every finite } \A_0 \in \A 
\mbox{ then } 
S\wedge\bigl(\bigwedge_{A\in \A} A \bigr) \not\leq \alpha^c.
$$

Further, since the preimage $F^{\gets}(y_0): X \to L$ of a point 
$y_0 \in Y$ under a fuzzy function $F: X \tto Y$ is obviously determined
by the equality 
$$F^{\gets}(y_0)(x) = \bigvee_{y\in Y} F(x,y)*y_0(y) = F(x,y_0),$$ 
the general definition of $(\alpha,\beta)$-compactness of an $L$-set in
this case means the following:

The preimage $F^{\gets}(y_0)$ of a point $y_0$ under a fuzzy function $F$ 
is $(\alpha,\beta)$-compact if for every family $\A$ of closed sets of
$X$ the following implication holds:
$$
\bigvee_x \Bigl(F(x,y_0)\wedge(\bigwedge_{A\in \A_0} A(x))\Bigr) 
\not\leq \beta^c\ \
\forall \A_0 \subset \A, |\A_0| < \aleph_0
$$
implies
$$\bigvee_x \bigl(F(x,y_0) \wedge
(\bigwedge_{A\in \A} A(x))\bigr) \not\leq \alpha^c.$$

Now we can introduce the following 

\begin{definition} 
A continuous fuzzy mapping $F\colon (X,E_X,\tau_X) \tto (Y,E_Y,\tau_Y)$ is
called $(\alpha,\beta)$-perfect if 
\begin{itemize}
\item 
$F$ is closed, i.e.\ $F^{\to}(A) \in \C_Y$ for every $A \in \C_X;$
\item 
the preimage $F^{\gets}(y)$ of every point $y \in Y$ is 
$(\alpha,\beta)$-compact.
\end{itemize}
\end{definition}

\begin{theorem} 
Let $F: (X,E_X,\tau_X) \tto (Y,E_Y,\tau_Y)$ be an 
$(\alpha,\gamma)$-perfect fuzzy function such that $\mu(F) = \top$ and 
$\sigma(F) = \top$. 
If the space $(Y,E_Y,\tau_Y)$ is $(\gamma,\beta)$-compact, then the space
$(X,E_X,\tau_X)$ is $(\alpha,\beta)$-compact.
\end{theorem}

\begin{proof}
Let $\A$ be a family of closed $L$-sets in $X$ such that 
$\bigwedge_{A\in\A_0} A \not \leq \beta^c$. 
Without loss of generality we may assume that $\A$ is closed under taking
finite meets. 
Let $B := B_A := F^{\to}(A)$ and let $\B := \{B_A: A \in \A \}$. 
Then, since $\mu(F) = \top$, by \ref{im-pr} (7) it follows that 
$B \not \leq \beta^c\ \forall B \in \B$, and moreover, since $\A$ is
assumed to be closed under finite meets,
$$
\begin{array}{lll}
B_{A_1}\wedge\ldots\wedge B_{A_n}&
=&
F^{\to}(A_1)\wedge\ldots\wedge F^{\to}(A_n)\\
&
\geq&
F^{\to}(A_1 \wedge\ldots\wedge A_n)\\
&
=&
F^{\to}(A),
\end{array}
$$
for some $A \in \A$, and hence 
$\bigwedge_{B\in\B_0}(B) \not \leq \beta^c$
for every finite subfamily $\B_0 \subset \B$. 
Hence, by $(\gamma,\beta)$-compactness of the space $(Y,E_Y,\tau_Y)$ we
conclude that $\bigwedge_{B \in \B}(B) \not \leq \gamma^c$, and therefore
there exists a point $y_0 \in Y$ such that 
$F^{\to}(A)(y_0) = B_A(y_0) \not \leq \gamma^c$ for all $A \in \A$. 
Now, applying $(\alpha,\gamma)$-compactness of the preimage
$F^{\gets}(y_0)$ and recalling that $\A$ was assumed to be closed under
taking finite meets, we conclude that 
$$\bigvee_x \bigl(F(x,y_0) \wedge 
(\bigwedge_{A\in \A} A(x))\bigr) 
\not\leq \alpha^c.$$
and hence, furthermore, 
$$\bigwedge_{A\in\A} A \not\leq \alpha^c.$$
\end{proof}

%\bibliographystyle{amsplain}
%\bibliography{26}
\providecommand{\bysame}{\leavevmode\hbox to3em{\hrulefill}\thinspace}
\providecommand{\MR}{\relax\ifhmode\unskip\space\fi MR }
\providecommand{\MRhref}[2]{%
  \href{http://www.ams.org/mathscinet-getitem?mr=#1}{#2}
}
\providecommand{\href}[2]{#2}

\end{document}